\documentclass[12pt]{amsart}
\usepackage{amssymb}
\theoremstyle{plain}

\newtheorem{thm}{Theorem}
\newtheorem{cor}[thm]{Corollary}
\newtheorem{lem}[thm]{Lemma}
\newtheorem{pro}[thm]{Problem}

\theoremstyle{definition}
\newtheorem{defi}[thm]{Definition}
\newtheorem{conj}[thm]{Conjecture}
\newtheorem{conv}[thm]{Convention}
\newtheorem{nota}[thm]{Notation}
\newtheorem{rem}[thm]{Remark}
\newtheorem{rems}[thm]{Remarks}
\newtheorem{exa}[thm]{Example}
\newtheorem{exas}[thm]{Examples}
\newtheorem{sit}[thm]{}

\newcommand{\brem}{\begin{rem}}
\newcommand{\brems}{\begin{rems}}
\newcommand{\erem}{\end{rem}}
\newcommand{\erems}{\end{rems}}
\newcommand{\bexa}{\begin{exa}}
\newcommand{\bexas}{\begin{exas}}
\newcommand{\eexa}{\end{exa}}
\newcommand{\eexas}{\end{exas}}
\newcommand{\bdefi}{\begin{defi}}
\newcommand{\edefi}{\end{defi}}
\newcommand{\bcor}{\begin{cor}}
\newcommand{\ecor}{\end{cor}}
\newcommand{\blem}{\begin{lem}}
\newcommand{\elem}{\end{lem}}
\newcommand{\bconv}{\begin{conv}}
\newcommand{\econv}{\end{conv}}
\newcommand{\bconj}{\begin{conj}}
\newcommand{\econj}{\end{conj}}
\newcommand{\bpro}{\begin{pro}}
\newcommand{\epro}{\end{pro}}
\newcommand{\bthm}{\begin{thm}}
\newcommand{\ethm}{\end{thm}}
\newcommand{\bnota}{\begin{nota}}
\newcommand{\enota}{\end{nota}}
\newcommand{\bsit}{\begin{sit}}
\newcommand{\esit}{\end{sit}}
\newcommand{\be}{\begin{eqnarray}}
\newcommand{\ee}{\end{eqnarray}}
\newcommand{\bproof}{\begin{proof}}
\newcommand{\eproof}{\end{proof}}

\def\ba{\begin{array}}
\def\ea{\end{array}}

\newcommand{\PP}{{\mathbb P}}
\newcommand{\R}{{\mathbb R}}
\newcommand{\C}{{\mathbb C}}
\newcommand{\Q}{{\mathbb Q}}

\title{Selected problems}

\author{Mikhail Zaidenberg}
\address{Universit\'e
Grenoble I, Institut Fourier, UMR 5582 CNRS-UJF, BP 74, 38402
St.\ Martin d'H\`eres c\'edex, France}
\email{zaidenbe@ujf-grenoble.fr}

\begin{document}



\maketitle

\bigskip

\noindent We start with the following complex-analytic version of
the famous Poincar\'e Conjecture. It was repeatedly  proposed in
\cite{Za1} and \cite{Za2}.

\noindent\bpro\label{pro1} Let $D$ be a strictly pseudoconvex
bounded domain in $\C^2$ with a smooth boundary $S$. Suppose $S$
is a homology sphere. Is it true that $S$ is diffeomorphic to
$S^3$?\epro

\noindent One may assume in addition that $S$ is real analytic or
even real algebraic, that $D$ is contractible, and ask whether
$D$ is a sublevel set of a real (strictly subharmonic) Morse
function with just one critical point. As a motivation, we
consider in {\em loc.cit.} an exhaustion of a smooth contractible
affine surface $X$ by a sequence of strictly pseudoconvex
subdomains $D_n\subseteq X$ with real-algebraic boundaries $S_n$.
If $X$ is not isomorphic to $\C^2$ (for instance, if $X$ is the
Ramanujam surface) then for all $n$ large enough, $S_n$ is a
homology sphere with a non-trivial fundamental group. Thus we
expect that such a domain $D_n$ cannot be biholomorphically
equivalent to a bounded domain in $\C^2$.

\medskip

The next 7 problems and conjectures also address contractible or
acyclic varieties; see \cite{Za2} for additional motivations.

\noindent\bpro\label{pro2} Consider a smooth contractible surface
$X$ of log-general type. Is the set of all algebraic curves in
$X$ with one place at infinity finite? Can it be non-empty?\epro

\noindent {\bf Rigidity Conjecture I} {\rm (Flenner-Zaidenberg
\cite{FlZa1})}. {\it Every smooth  $\Q$-acyclic surface  $X$  of
log-general type is rigid and has unobstructed deformations. In
particular, for any minimal smooth completion $V$ of $X$ by a
divisor $D$ with simple normal crossings, $\chi(\Theta_{V\langle
D\rangle})=K_V (K_V + D) = 0$, or equivalently\footnote{As was
conjectured by T. tom Dieck \cite{tD}.} $3{\overline c}_2 (X) -
{\overline c}_1^2 (X) = 5$.}

\medskip

\noindent This conjecture has been verified in a number of cases
\cite{FlZa1}. Of special interest are the $\Q$-acyclic surfaces of
type $X=\PP^2\setminus C$, where $C$ is a rational cuspidal plane
curve\footnote{By a {\it cusp}  of $C$ we mean any locally
irreducible (i.e., unibranch) singularity of $C$.} (see e.g.
\cite{FlZa2, FlZa3, OrZa}).

\medskip

\noindent {\bf Rigidity Conjecture II} {\rm (Flenner-Zaidenberg
\cite{FlZa3})}.  {\it Every rational cuspidal plane curve with at
least 3 cusps is projectively rigid and has unobstructed
deformations. Thus given $d$, there is, up to projective linear
transformations, at most countable set of such curves of degree
$d$. If $m$ stands for the maximal multiplicity of singular
points of $C$ then, moreover, $m\ge d-4$, besides a finite number
of possible exceptions.}

\medskip

\noindent Presumably, a rational cuspidal plane curve $C$ cannot
have more than 4 cusps (cf. \cite{Or2}). All such curves with
$m\ge d-3$ were classified in \cite{FlZa2, FlZa3, SaTo}, and in
\cite{Fe} all those with $m= d-4$ and with unobstructed
deformations. See also the very interesting recent preprint
\cite{FLMN}.

\medskip

\noindent {\bf Finiteness Conjecture} \cite{Za2}. {\it For a
$\Q$-acyclic surface $X$ of log--general type, we consider a
minimal smooth SNC--completion $V$ of $X$ with the boundary
divisor $D$. We conjecture that the set of all possible
Eisenbud--Neumann diagrams\footnote{We recall that $\breve
\Gamma_D$ is a tree without vertices of valence 2, obtained from
the dual graph $\Gamma_D$ of $D$ by contracting all its linear
chains.} $\breve \Gamma_D$ of such curves $D$ is finite. }

\medskip

\noindent See \cite{Or1} for a result confirming the conjecture.

\medskip

\noindent  {\bf Isotriviality Conjecture} \cite{Za2}. {\it
Consider a (reduced, irreducible) hypersurface $X = \{p = 0 \}$
in $\C^n$, where $p \in \C[x_1 ,\dots, x_n]$. If $X$ is
contractible then $p$ is an isotrivial polynomial, i.e. the
generic fibres of $p$ are pairwise isomorphic. }

\medskip

\noindent Does there exist any smooth such hypersurface $X$ of
log-general type?

\medskip

\noindent  {\bf Quasihomogeneity Conjecture} \cite{Za2}. {\it If
$X$ as above has an isolated singular point then, up to a
polynomial automorphism of $\C^n$, $p$ is equivalent to a
quasihomogeneous polynomial.}

\medskip

By an exotic $\C^n$  we mean a smooth affine algebraic variety
diffeomorphic to $\R^{2n}$ and non-isomorphic to $\C^n$ (see
\cite{Za3}).

\noindent \bpro\label{pro8} Does there exist a pair $(X,Y)$ of
exotic $\C^n$'s which are biholomorphic and non-isomorphic? What
can be said if $X=\C^n$? Does there exist a non-trivial
deformation family of exotic $\C^n$'s with the same underlying
analytic structure? Show that the Russell cubic\footnote{Which is
known to be an exotic $\C^3$ \cite{ML}.} $X = \{x + x^2y + z^3 +
t^2 = 0 \} \subset \C^4$ is not biholomorphic to $\C^3$. Does
there exist any non-rational exotic $\C^n$?\epro

\noindent The latter question was proposed by Hirzebruch and Van
de Ven \cite{vDV}; the negative answer is known in dimension 2
(Gurjar-Shastri; see \cite{GPS} for a more general result).

\medskip

We turn further to the Abhyankar-Sathaye Problem on equivalence of
regular embeddings $\C^k\hookrightarrow \C^n$. Although many
conjectural counterexamples were proposed (see e.g. \cite{As,
KaZa1, KaVeZa}), so far the problem resists due to the lack of
suitable invariants. A particular class of such examples, called
the V\'en\'ereau polynomials (see e.g. \cite{KaZa4}), leads to the
following related problem.

\noindent \bpro\label{pro9} Classify the algebraic fiber bundles
over the punctured plane $\C^2\setminus\{\bar 0\}$ with vector
space fibers. \epro

\noindent  Let $\Upsilon$ be such a vector bundle with $\C^n$
($n>1$) as the standard fiber. Then the structure group, say, $G$
of $\Upsilon$ is a subgroup of the infinite dimensional group of
polynomial automorphisms of $\C^n$. In case that $G$ reduces to
the subgroup of affine transformations, by a result of Brenner
\cite{Br} the bundle $\Upsilon$ is trivial if and only if its
total space is not affine. We wonder (after H. Brenner) whether
this criterion still holds in the general case.

\noindent \bpro\label{pro10} {\rm (Kaliman-Zaidenberg
\cite{KaZa3})} Let $f:X\to S$ be a quasiprojective flat family
with generic fiber $\C^n$. Show that $X$ contains a cylinder,
i.e. that $f$ is equivalent to a trivial family $U\times\C^n\to
U$ over a Zariski open subset $U\subseteq S$. \epro

\noindent This is true for $n\le 2$ (see e.g. \cite{KaZa3}). The
problem is a weakened form of the Dolgachev-Weisfeiler
Conjecture, which suggests that $(X,f,S)$ should be a fiber
bundle assuming $S$ smooth and all the fibers of $f$ being
reduced and isomorphic to $\C^n$.

\medskip

Denote by $\delta_n(X,p)$ the Watanabe plurigenera of an isolated
normal singularity $(X,p)$ of a complex space (see e.g.,
\cite{FlZa3}).

\noindent \bpro\label{pro11} {\rm (Flenner-Zaidenberg
\cite{FlZa3})} Let $X$ be a normal quasiprojective variety with an
isolated singular point $p$. Is it true that, whenever $X$ admits
a non-trivial regular $\C_+$-action, then $\delta_n(X,p)=0$ for
all $n$ and, moreover, $(X,p)$ is a quotient singularity? Show
that the 3-fold $x^3+y^3+z^3+t^3=0$ in $\C^4$ does not admit such
an action. Describe all (quasi)homogeneous hypersurfaces $X$ in
$\C^n$ with a $\C_+$-action, and, more generally, those $X$ which
admit a non-constant polynomial map $\C\to X$ whose image is not
an orbit closure of a $\C^*$-action on $X$. \epro

\noindent \bpro\label{pro12} {\rm (Flenner-Zaidenberg
\cite{FlZa5})} Find an example of an affine surface $X$ with a
trivial Makar-Limanov invariant\footnote{That is $X$ admits two
non-trivial $\C_+$-actions with transversal general orbits.},
which does not admit any non-trivial $\C^*$-action. Classify all
the normal affine 3-folds (respectively, $n$-folds) with (a torus
action and) a trivial Makar-Limanov invariant. \epro

\noindent All normal affine surfaces with a  $\C^*$-action and a
trivial Makar-Limanov invariant were described in \cite{FlZa5}.
As for the tori actions on affine varieties, see e.g. \cite{AlHa};
cf. also \cite{KaZa2}.

\medskip

In another direction, we address the following

\noindent \bpro\label{pro123} {\rm (Shiffman-Zaidenberg
\cite{ShZa1, ShZa2})} Construct an example of a Kobayashi
hyperbolic projective surface in $\PP^3$ of degree 5. Construct,
for every $n>2$, an example of a Kobayashi hyperbolic
hypersurface in $\PP^n$ of degree $d_n$ with $d_n\le Cn$, where
$C>0$ does not depend on $n$ (respectively, with $d_n=2n-1$).\epro

\noindent An example of a hyperbolic surface in $\PP^3$ of degree
6 was recently found by Duval \cite{Du}, whereas the least
admissible degree $d$ of such a surface  is 5. In examples of
hyperbolic hypersurfaces in higher dimensions, the best known
asymptotics for the degrees so far is $d\sim n^2$, see e.g.
\cite{ShZa1, Za4}.

\end{document}